\documentclass[a4paper,twoside,12pt]{article} 
\usepackage{amsfonts,amssymb,mathrsfs,amsmath}
\usepackage[dvips]{color}
\usepackage{graphicx}

\title{A Tiling Proof of Binomial Identities related to the Lucas cube}
\date{}
\author{Andrzej P. Kisielewicz\\
\\
{\small Wydzia{\l} Matematyki, Informatyki i Ekonometrii, Uniwersytet Zielonog\'orski}\\
{\small ul. Z. Szafrana 4a, 65-516 Zielona G\'ora, Poland}\\
{\small A.Kisielewicz@wmie.uz.zgora.pl}\\
}
\vspace{-5 mm}
\numberwithin{equation}{section}

\newtheorem{df}{\sc Definition}

\newtheorem{uw}{\sc Remark}
\newtheorem{uwi}[uw]{\sc Remarks}

\newtheorem{nap}{\sc Example }

\newtheorem{nps}[nap]{\sc Examples}

\def\ka #1{\mathscr{#1}}
\def\kal #1 #2{\mathscr{#1}^{#2}}

\def\zet{\mathbb{Z}}

\def\er{\mathbb{R}}

\begin{document}

\numberwithin{pr}{section}
\numberwithin{uw}{section}
\maketitle
\begin{abstract}
Using a cube tiling of $\er^n$ constructed by Lagarias and Shor  a tiling proof of three well-known  binomial identities related to the Lucas cube is given. 


\end{abstract}

\section{Introduction}

%\maketitle

A cube tiling of $\er^n$ is a family of cubes $[0,2)^n+T=\{[0,2)^n+t:t\in T\}$, where $T\subset \er^n$, which fill in the whole space without gaps and overlaps. 
%Every partition code $V\subset S^d$ can be used to construct a cube tiling of $\er^d$ (see e.g. ?).    
In \cite{LS2} Lagarias and Shor constructed a cube tiling code of $\er^n$. In this note we show how it can be used to prove the following well-known identities:
%on the other two constructions of cube tilings of $\er^d$.
%which can by used in many areas of the tiling theory. 
\begin{equation}
\label{44}
\sum_{k\geq 0}{{n-k}\choose {k}} \frac{n}{n-k}2^k=2^n+(-1)^n, 
\end{equation}
\begin{equation}
\label{45}
\sum_{k\geq 0}{{n-k}\choose {k}}2^k=1/3(2^{n+1}+(-1)^{n}) 
\end{equation}
and
\begin{equation}
\label{46}
\sum_{k\geq 0}{{n-k}\choose {k}} \frac{k}{n-k}2^k=1/3(2^n+(-1)^n2).  
\end{equation}

%In most cases we have kept the original notation of the paper \cite{LS2}. In particular, we denote the dimension by  $n$ instead of $d$. 
The code of Lagarias and Shor is constructed as follows.
Let $n\geq 3$ be an odd positive integer, and let $A$ be an $n\times n$ circulant matrix of the form $A=A(n)={\rm circ}(1,2,0,\ldots,0)$. Let $A^T$ be the transpose of $A$. By $V(A)$ and $V(A^T)$ we denote the sets of all the sums of the different rows in $A$ and $A^T$, respectively. Moreover, we add to these sets the vector $(0,...,0)$.  Let
$$
V=V_e(A)\cup (V_o(A^T)+(2,\ldots ,2)) \mod 4,
$$
where  $V_e(A)$ denotes the set of all vectors in $V(A)$ with an even number of $3$s,  and $V_o(A^T)$ is the set of all vectors in $V(A^T)$ with an odd number of $0$s.
We will refer to the code $V$ as the {\it Lagarias-Shor cube tiling code}. This code has very interesting applications. Originally in \cite{LS2} it was used to design a certain cube tiling of $\er^n$ that was the basis for estimating distances between cubes in cube tilings of $\er^n$. Recently in \cite{K} the Lagarias-Shor cube tiling code was used to  construct interesting partitions and matchings of an $n$-dimensional cube. %The same part can be also used to construct an interesting systems of abstract boxes (so-called {\it rigid plyboxes}). 

To obtain a cube tiling of $\er^n$ from the code $V$, let $T=V-{\bf 1}+4\zet^n$, where ${\bf 1}=(1,\ldots ,1)$.  
It follows from  \cite[Proposition 3.1 and Theorem 4.1]{LS2} that $[0,2)^n+T$ is a cube tiling of $\er^n$. (To be precisely, a tiling considered in \cite{LS2} is of the form $[0,1)^n+T'$, where $T'=\frac{1}{2}V+2\zet^n$, but  $[0,2)^n+T=[0,2)^n+2T'-{\bf 1}$.)

All three identities (\ref{44})--(\ref{46}) are related to the {\it Lucas cube $\Lambda_n$}. This is a graph whose vertices are all elements of  $\{0,1\}^n$ which do not contain two consecutive 1s as well as vertices having $1$ at the first and last position simultaneously.  It is known that ${{n-k}\choose {k}} \frac{n}{n-k}$ is the number of all vertices $v$ in the  Lucas cube $\Lambda_n$ of weight $k$, i.e., containing $k$  1s. This is also the number of all $k$-element subsets of the set $[n]=\{1, ...,n\}$ without two consecutive integers and which do not contain the pair $1$, $n$ (\cite{Mu}). The number  of all vertices in $\Lambda_n$ of weight $k$ which have $1$ at the $i$th position, $i\in [n]$, is equal to ${{n-k}\choose {k}} \frac{k}{n-k}$, while ${{n-k}\choose {k}}$ is the number of all  vertices in $\Lambda_n$ of weight $k$ which have $0$ at the $i$th position.

%The identity (\ref{44}) is closely related to the {\it Lucas cube $\Lambda_n$}. This is a graph whose vertices are all elements of  $\{0,1\}^n$ which do not contain two consecutive 1s as well as vertices having $1$ at the first and last position simultaneously.  It is known that ${{n-k}\choose {k}} \frac{n}{n-k}$ is the number of all vertices $v$ in the  Lucas cube $\Lambda_n$ of weight $k$, i.e., containing $k$  1s. This is also the number of all $k$-element subsets of the set $[n]=\{1, ...,n\}$ without two consecutive integers and which do not contain the pair $1$, $n$ (\cite{Mu}).  

%This is a graph whose vertices are all elements of  $\{0,1\}^n$ which do not contain two consecutive 1s as well as vertices having $1$ at the first and last position simultaneously.  
%(Clearly, from $U(n)$ we can obtain the set of vertices in $L_n$ if we make the substitutions: $*\rightarrow 1$, $2\rightarrow 0$ and $0\rightarrow 0$.) 
%Similarly, vertices in the {\it Fibonacci cube $\Gamma_n$} are all elements of  $\{0,1\}^n$ which do not contain two consecutive 1s  
%allowing this time vertices with two 1s at the first and last position. The number of all vertices of weight $k$ in the Fibonacci cube $\Gamma_n$ is equal to ${{n-k+1}\choose {k}}$ (\cite{Mu1}). 

\section{Tiling proofs}

Since the Lagarias-Shor cube tiling code is defined for odd numbers, we prove identities (\ref{44})-(\ref{46}) for odd and even positive integers separately. %We start with the proof of 

\smallskip
\noindent
\textit{Proof of (\ref{44}) for $n\geq 3$ odd.} We intersect the cube $[0,2)^n$ with the cubes from the tiling  $[0,2)^n+T$. Let $\ka F(n)=\ka F=\{[0,2)^n\cap ([0,2)^n+t): t\in T\}$. Since $[0,2)^n+T$ is a tiling, $\ka F$ is a partition of the cube $[0,2)^n$ into boxes. Let $m(K)$ denote the volume of the box $K\in \ka F$, and let $m(\ka F)=\sum_{K\in \ka F}m(K).$ For every $K\in \ka F$ we have $m(K)=2^k$, where $k$ is the number of 1s in the vector $v\in V$ such that $K=[0,2)^n\cap ([0,2)^n+v-{\bf 1})$. Let $M_k=|\{K\in \ka F: m(K)=2^k\}|$. The family $\ka F$ is a partition of $[0,2)^n$ and therefore $m(\ka F)=2^n$ and 
\[
m(\ka F)-2=\sum_{k\geq 1}M_k2^k.
\]
Note now that if $v\in V$ contains 3 at some position $i\in [n]$, then the cubes $[0,2)^n+v-{\bf 1}$ and $[0,2)^n$ are disjoint. This means that these two cubes intersect if and only if $v\in U\cup \{(0,\ldots ,0), (2,\ldots ,2)\}$,  where $U$ consists of all sums of non-adjacent rows of the matrix $A$ where the first and last rows are treated as adjacent.  Thus, $M_k$ is the number of all $k$-element subsets of the set $\{1, ...,n\}$ without two consecutive integers and which do not contain the pair $1$, $n$. Hence, $M_k={{n-k}\choose {k}} \frac{n}{n-k}$. This completes the proof of $(\ref{44})$ for $n\geq 3$ odd. 
\hfill{$\square$}
%(Observe that, making the substitution $2\rightarrow 0$ in every $v\in  U\cup \{(0,\ldots ,0), (2,\ldots ,2)\}$ we get the set of vertices of the Lucas cube $\Lambda_n$.)  

This proofs needs only the portion $U=U(n)$ of the Lagarias-Shor cube tiling code, where the code $U$ is as in the previous proof: it consists of all sums of non-adjacent rows of the matrix $A(n)$, where the first and last rows are treated as adjacent. For $n=5$ we have   

\medskip
\begin{displaymath}
A(5)=
\left[\begin{array}{ccccc}
1 & 2 & 0 & 0  & 0\\
0 &1 &2 & 0 &0\\
0 &0 &1 &2 &0\\
0 &0 &0 &1 &2\\
2 &0 &0 &0 &1\\
\end{array}\right],\;\;\;  
U(5)=
\left[\begin{array}{ccccc}
1 &2 &0 &0 &0\\
0 &1 &2 &0 &0\\
0 &0 &1 &2 &0\\
0 &0 &0 &1 &2\\
2 &0 &0 &0 &1\\
1 &2 &1 &2 &0\\
1 &2 &0 &1 &2\\
0 &1 &2 &1 &2\\
2 &1 &2 &0 &1\\
2 &0 &1 &2 &1\\
\end{array}\right].
\end{displaymath}

\medskip
This is easily seen that making in every vector $v\in  U\cup \{(0,\ldots ,0)\}$ the substitution $2\rightarrow 0$ we obtain the set of all vertices in the Lucas cube $\Lambda_n$. 

\medskip
To  prove (\ref{45}) and (\ref{46}) for $n\geq 3$ odd let $\ka F^i_0, \ka F^i_2$ and $\ka F^i_1$,  $i\in [n]$ denote the sets of all boxes in $\ka F$ which have the factors $[0,1),[1,2)$ and $[0,2)$ at the $i$th position, respectively. Since $\ka F=\{[0,2)^n\cap ([0,2)^n+v-{\bf 1}): v\in U\cup \{(0,\ldots ,0), (2,\ldots ,2)\} \}$,
%every box in $\ka F$ is determined by $v\in U\cup \{(0,\ldots ,0),(2,\ldots ,2)\}$. Thus, 
for every $k\in \{0,1,2\}$ the set $\ka F^i_k$ consists of all boxes in $\ka F$ which are determined by the vectors $v\in U\cup \{(0,\ldots ,0),(2,\ldots ,2)\}$ such that $v_i=k$.
Let 
$$
m(\ka F^i_{02})=\sum_{K\in \ka F^i_{02}}m(K)\;\; {\rm and} \;\; m(\ka F^i_1)=\sum_{K\in \ka F^i_{1}}m(K),
$$
where $ \ka F^i_{02}=\ka F^i_{0}\cup \ka F^i_{2}$.

The partition $\ka F$ (Figure 1) has the structure which will be utilized below. Note that for every $i\in [n]$ the set  $\bigcup \ka F^i_{02}$,  in an {\it $i$-cylinder}, i.e., for every line segment $l_i\subset [0,2)^n$ of length 2 which is parallel to the $i$th edge of the cube $[0,2)^n$ we have 
\begin{equation}
\label{wal}
l_i\subset \bigcup \ka F^i_{02}\;\; {\rm or} \;\; l_i\cap \bigcup \ka F^i_{02}=\emptyset.
\end{equation}
Obviously, the set $\bigcup \ka F^i_{1}$ is an $i$-cylinder too.

\medskip
{\center
\includegraphics[width=10cm]{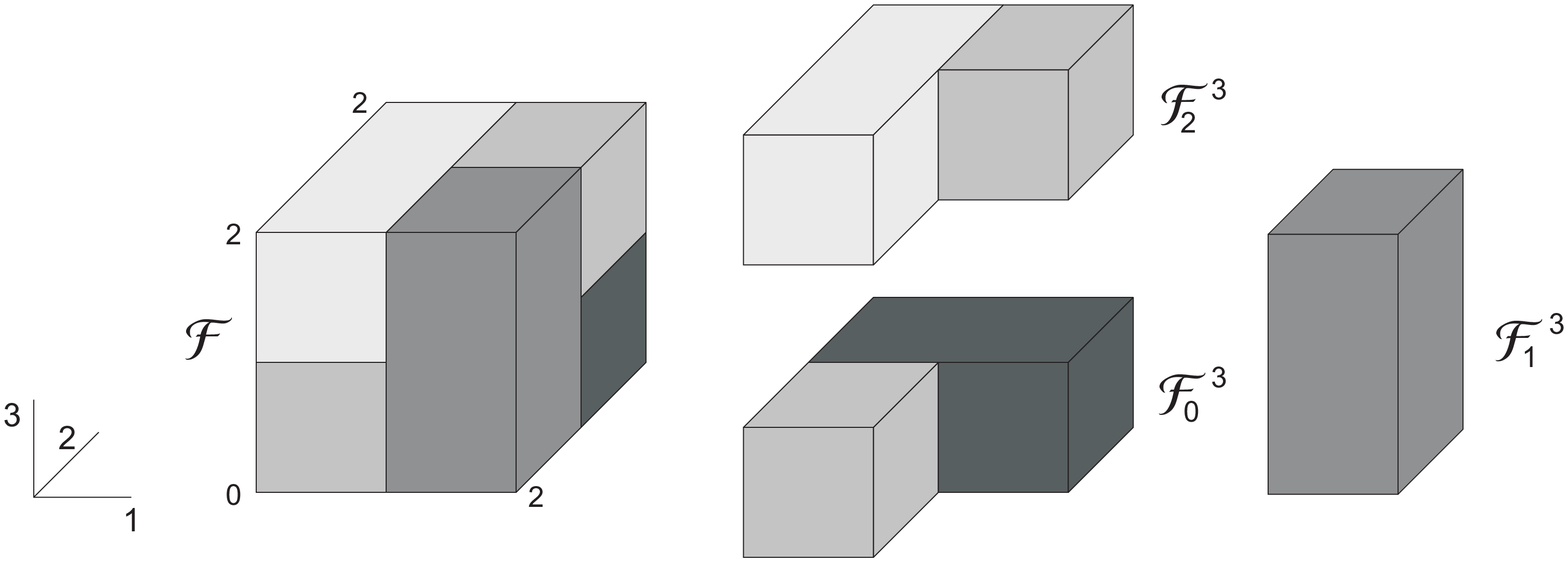}\\
}

\bigskip
\noindent{\footnotesize Fig. 1. %The circulant matrix $A(3)={\rm circ}(1,2,0)$ and the partition $\ka F$ of $[0,2)^3$. 
The boxes in $\ka F$ are determined by the vectors $U=\{(1,2,0),(0,1,2),(2,0,1)\}$ (the three "long" boxes) and $\{(0,0,0), (2,2,2)\}$ (the two unit cubes). 

}

\medskip
\noindent
\textit{Proof of (\ref{45}) and (\ref{46}) for $n\geq 3$ odd.} We will calculate $m(\ka F^i_{02})$ and $m(\ka F^i_1)$. For every $v\in U$ we have $v_i=1$ if and only if $v_{i+1}=2$. Thus,  $\ka F^{i+1}_{2}=\ka F^i_1\cup \{[1,2)^n\}$ and then $m(\ka F^{i+1}_{2})=m(\ka F^i_1)+1$ (clearly, $n+1$ is taken modulo $n$). It follows from (\ref{wal}) that $m(\ka F^{i}_0)=m(\ka F^{i}_2)$ (see Figure 1). 
%Any realization of the set of words $U^{i+1}_{02}$ is an $(i+1)$-cylinder. This implies that $U^{i+1}_0$ and $U^{i+1}_2$ are equivalent codes, and thus $m(U^{i+1}_0)=m(U^{i+1}_2)$. 
Since $A$ is a circulant matrix, we have $m(\ka F^i_1)=m(\ka F^{j}_1)$ for $i,j\in [n]$, and $m(\ka F)=m(\ka F^i_0)+m(\ka F^i_2)+m(\ka F^i_1)$ because $\ka F$ is a partition. Thus, $2^n=3m(\ka F^i_1)+2$ and consequently 
\begin{equation}
\label{uu} 
 m(\ka F^i_{02})=2/3(2^n+1)\;\;\ {\rm and}\;\;\; m(\ka F^i_1)=1/3(2^n-2).  
\end{equation}
As it was noted before the proof, the code $U\cup \{(0,\ldots ,0)\}$ can be identify with the set of all vertices in the Lucas cube $\Lambda_n$. %Thus, in what follows we can use the enumerative properties of $\Lambda_n$. 

%Let $\Lambda_n^{1,k}$ be the set  of all vertices of weight $k$ in the Lucas cube $\Lambda_n$ which have $1$ at the first position. 
%We have  $|\Lambda_n^{1,k}|=|\ka F^1_1|$.
%We will calculate $|L_n^{1,k}|$. 
%Observe that, $(1,w_2,\ldots,w_n)\in \Lambda_n^{1,k}$ if and only if $(w_3,\ldots ,w_{n-1})$ is a vertex of the Fibonacci cube $\Gamma_{n-3}$ of weight $k-1$. 
%Thus, $|\Lambda_n^{1,k}|={{n-k-1}\choose {k-1}}.$ Since $\frac{n-k}{k}{{n-k-1}\choose {k-1}}={{n-k}\choose {k}}$, 
Since ${{n-k}\choose {k}}\frac{k}{n-k}$ is the number of all vertices of weight $k$ in $\Lambda_n$ with $1$ at the first position, we have 
%Recall that for every $i\in [n]$ the set $U^i_*$ consists of words having the star at the $i$th position. Thus, by the above equality, we have
$|\ka F_1^1|= \sum_{k\geq 1}{{n-k}\choose {k}}\frac{k}{n-k}$, and consequently
%$|\ka F^1_1|=\sum_{k=1}^n{{n-k}\choose {k}}\frac{k}{n-k}$ and consequently 
$$
m(\ka F^1_1)=\sum_{k\geq 1}{{n-k}\choose {k}}\frac{k}{n-k}2^k,
$$
which, by (\ref{uu}), gives (\ref{46}) for $n\geq 3$ odd. 
%Similarly, by the first formula in (\ref{vv}), we obtain (\ref{46}) for $n-1\geq 2$ even. 
Since 
$$
\sum_{k\geq 0}{{n-k}\choose {k}} \frac{n}{n-k}2^k=\sum_{k\geq 0}{{n-k}\choose {k}}2^k +\sum_{k\geq 0}^{n}{{n-k}\choose {k}} \frac{k}{n-k}2^k,
$$
%\begin{equation}
%\label{inter}
%$$
%\sum_{k=0}^{n}{{n-k}\choose {k}} \frac{n}{n-k}2^k=\sum_{k=0}^{n}{{n-k}\choose {k}}2^k +\sum_{k=0}^{n}{{n-k}\choose {k}} \frac{k}{n-k}2^k,
%$$
%\end{equation}
the proof of the identity (\ref{45}) for $n\geq 3$ odd is also completed. 
\hfill{$\square$}

\medskip 
%These identities are close related to the partition $\ka U$.
For $n\geq 3$ odd all three identities are strongly related to the partition $\ka F$. The sums $\sum_{k\geq 1}{{n-k}\choose {k}}2^k+2$ and $\sum_{k\geq 1}{{n-k}\choose {k}}\frac{k}{n-k}2^k$ are the total volumes of the boxes from the partition $\ka F$ which belong to the sets $\ka F^i_{02}$ and $\ka F^i_1$, respectively (the number $2$ in the first sum is the sum of the volumes of the boxes $[0,1)^n$ and $[1,2)^n$). The summands ${{n-k}\choose {k}}2^k$ and ${{n-k}\choose {k}}\frac{k}{n-k}2^k$ for  $k=1,\ldots ,{\left\lfloor \frac{n}{2}\right\rfloor}$ are the total volumes of the boxes in $\ka F^i_{02}$ and $\ka F^i_1$, respectively which have exactly $k$ factors $[0,2)$. %It is interesting that the natural representation
%Now this is seen that the identities (\ref{45}) and (\ref{46}) are also related to the Lucas cube $\Lambda_n$: the number $\sum_{k=1}^{n}{{n-k}\choose {k}}$ counts all vertices in $\Lambda_n$ with $0$ at the fixed position $i\in [n]$, and 

%$$
%\sum_{k=1}^{n}{{n-k}\choose {k}} \frac{n}{n-k}2^k+2=(\sum_{k=1}^{n}{{n-k}\choose {k}}2^k+2) +\sum_{k=1}^{n}{{n-k}\choose {k}} \frac{k}{n-k}2^k
%$$
%coincides with the natural representation of the partition $\ka F$ into two $i$-cylinders: $\cup \ka F=\cup \ka F^i_{02}\cup \cup \ka F^i_1$ (see Figure 1).  

\bigskip
The identities (\ref{44})--(\ref{46}) for $n-1\geq 2$ even can be derived from the partitions $\ka F(n)$ and $\ka F(n-2)$, where $\ka F(1)=\{[0,1),[1,2)\}$.

\medskip
\noindent
\textit{Proofs of (\ref{44})-(\ref{46}) for  $n-1\geq 2$ even.}
Denote by $r_1,\ldots, r_n$ the rows of the matrix $A(n)$, and let  $\ka G=\ka G(n-1)\subset \ka F(n)$ be the set  of all boxes which are determined by the vectors $v\in U(n)$ which are sums of non-adjacent rows from the set $\{r_1,\ldots ,r_{n-1}\}$, where  $r_1$ and $r_{n-1}$ are treated as adjacent.
%except for all words starting with the string $12$ and ending with $12$. This means that these sums do not contain the rows $r_1$ and $r_{n-1}$ simultaneously. 
Thus, the number $|\ka G|$ is the same as the number of all vertices in the Lucas cube $\Lambda_{n-1}$. Consequently 
$$
m(\ka G)=\sum_{k\geq 1}{{n-1-k}\choose {k}} \frac{n-1}{n-1-k}2^k.
$$ 
Since $\ka G^n_0=\ka F^n_0\setminus \{[0,1)^n\}$, it follows that $m(\ka G^n_0)=m(\ka F^n_0)-1$, and by (\ref{uu}) and the fact that $m(\ka F^n_0)=m(\ka F^n_2)$ we get 
$$
m(\ka G^n_0)=2/3(2^{n-1}-1).
$$
We now calculate  $m(\ka G^n_2)$. Every box in $\ka G^n_2$ is generated by a vector $v\in U$ which has $2$ at the $n$th position. Therefore, $v=r_{n-1}+\sum_{i\in I} r_i$ for some $I\subset \{2,\ldots, n-3\}$.
%To determine $m(V^n_2)$ notice that every word in the set $V^n_2\setminus \{r_{n-1}\}$ is of the form $r_{n-1}+\sum_{i\in I} r_i$ for some $I\subset \{2,\ldots, n-3\}$. 
Let $R$ be the set of all such sums $\sum_{i\in I} r_i$. Every word in $R$ is a sum of non-adjacent rows from the set $\{r_2,\ldots ,r_{n-3}\}$, where $r_2$ and $r_{n-3}$ are not treated as adjacent. Let $U^{n-2}_0(n-2)$ be the set of all vectors in $U(n-2)$ having 0 at the last position. Observe now that the function $b:R\rightarrow U^{n-2}_{0}(n-2)$ defined by the formula $b(u)=\sum_{i\in I-1}h_i$, where $h_1,\ldots, h_{n-2}$ are rows of the matrix $A(n-2)$ and $I-1=\{i-1:i\in I\}$, is a bijection. 
%Moreover, $\bar{g}(s_*,u)=\bar{g}(p_*,b(u))$, where the word $p_*$ consists of $n-2$ stars. 
Therefore, $m(\ka G^n_2)=2m(\ka F^{n-2}_0(n-2))$ (recall that we add $r_{n-1}$ to $\sum_{i\in I} r_i$).  By (\ref{wal}), $m(\ka F^{n-2}_0(n-2))=m(\ka F^{n-2}_2(n-2))$, and 
% that is $m(R)=1/3(2^{n-2}+1)-1$. By the relationship between the sets $V^n_2$ and $R$, we have $m(V^n_2)=2m(R)+2$, and hence $ m(V^n_2)=1/3(2^{n-1}+2)$.
by (\ref{uu}),  
$$
m(\ka G^n_2)=1/3(2^{n-1}+2).
$$
Thus, $m(\ka G)=m(\ka G^n_0)+m(\ka G^n_2)=2^{n-1}$ because $\ka G^n_1=\emptyset$. This completes the proof of (\ref{44}) for $n-1\geq 2$ even. 

Since $m(\ka G^{i+1}_2)=m(\ka G^i_1)$, $m(\ka G^{i}_1)=m(\ka G^j_1)$ and $m(\ka G)=m(\ka G^i_1)+ m(\ka G^i_{02})$ for $i,j\in [n-1]$, it follows that 
$$
m(\ka G^i_1)=1/3(2^{n-1}+2)\;\; {\rm and}\;\;  m(\ka G^i_{02})=2/3(2^{n-1}-1)
$$
 for $i\in [n-1]$. 
 
By the definition of the set $\ka G$, we have $|\ka G^1_1|=\sum_{k\geq 1}{{n-1-k}\choose {k}}\frac{k}{n-1-k}$, and thus
$$
m(\ka G^1_1)=\sum_{k\geq 1}{{n-1-k}\choose {k}}\frac{k}{n-1-k}2^k 
$$
which proves (\ref{46}) for $n-1\geq 2$ even. Having this, in the same manner as for $n\geq 3$ odd we prove $(\ref{45})$ for $n-1\geq 2$ even.
\hfill{$\square$}

\begin{uw}{\rm
There are many tiling proofs that rely on counting the number of $1$-dimensional tilings of a $1\times n$ board  by  polyominoes (squares, dominoes etc.) (see \cite{BQ,P2}). In our case we examine just one $n$-dimensional tiling of $\er^n$ by translates of the cube $[0,2)^n$, and especially it is exploited the structure of that tiling.}
\end{uw} 
%At the end let us mentioned that in many tiling proofs we 
%The method, which is used in many proofs in \cite{P1,BQ,P2}, relies on counting the number of $1$-dimensional tilings by squares and dominoes.
%In \cite{P1,BQ,P2} proofs are relied on counting of the number of $1$-dimensional tilings by squares and dominos. 
%In our case, we examine only one $n$-dimensional tiling (a code of a cube tiling of $\er^n$, to be precise), and   especially it is extremely exploited the way of generating of it. 
At the end we show that for $n\geq 3$ odd the set of all vertices of the Lucas cube $\Lambda_n$ is a selector of a discrete analogue of the partition $\ka F$ 

\section{Vertices of the Lucas cube as a selector}

Let $L\cup \{(0,\ldots ,0), (1,\ldots ,1)\}$ be the code that arises from $U$ by making in every vector $v\in  U\cup \{(0,\ldots ,0), (2,\ldots ,2)\}$ the following substitutions:  $0\rightarrow 0$, $2\rightarrow 1$ and $1\rightarrow *$
For example, for $n=5$ we have %we get the set of vertices of the Lucas cube $\Lambda_n$.

\medskip
\begin{displaymath}
L=
\left[\begin{array}{ccccc}
* &1 &0 &0 &0\\
0 &* &1 &0 &0\\
0 &0 &* &1 &0\\
0 &0 &0 &* &1\\
1 &0 &0 &0 &*\\
* &1 &* &1 &0\\
* &1 &0 &* &1\\
0 &* &1 &* &1\\
1 &* &1 &0 &*\\
1 &0 &* &1 &*\\
\end{array}\right].
\end{displaymath} 

\medskip
The set $L$ consists of all sums of non-adjacent rows of the matrix ${\rm circ}(*,1,0,\\ \ldots ,0)$, where the first and last rows are treated as adjacent. Therefore, making the substitution $*\rightarrow 0$ in every vector of  $L\cup \{(0,\ldots ,0)\}$ we obtain the set $V(\Lambda_n)$ of all vertices in the Lucas cube $\Lambda_n$.

The code $L\cup \{(0,\ldots ,0),(1,\ldots ,1)\}$ induces a partition $\ka L$ of discrete box $\{0,1\}^n$ into boxes which is a discrete analogue of the partition $\ka F$ from the previous section. The boxes $K(l)=K_1(l)\times \cdots \times K_n(l)\in \ka L$, where $l\in  L\cup \{(0,\ldots ,0),(1,\ldots ,1)\}$, are of the form:  

$$
K_i(l)=\left\{ \begin{array}{rl}
\{0\} & \textrm{if $l_i=0$,} \\
\{1\} & \textrm{if $l_i=1$,}\\
\{0,1\}  & \textrm{if $l_i=*$}\\
\end{array} \right.
$$
for $i\in [n]$.
  
Observe now that for every $n\geq 3$ odd the set $V(\Lambda_n)$ of the vertices of the Lucas cube is a {\it selector} of the family of boxes $\ka L\setminus \{\{1\}\times \cdots \times \{1\}\}$: for every $v\in V(\Lambda_n)$ there is exactly one $K(l)\in  \ka L\setminus \{\{1\}\times \cdots \times \{1\}\}$
such that 
$$
v\in K(l).
$$
Indeed, let $K(l)\in \ka L\setminus \{\{1\}\times \cdots \times \{1\}\}$ and pick a point $v=v_i\cdots v_d\in K(l)$ in the following way:
$$
v_i=\left\{ \begin{array}{rl}
0 & \textrm{if $K_i(l)=\{0\}$,} \\
1 & \textrm{if $K_i(l)=\{1\}$,}\\
0  & \textrm{if $K_i(l)=\{0,1\}$.}\\
\end{array} \right.
$$
 Since $l$ does not contain two consecutive 1s modulo n and  if $K_i=\{0,1\}$, then $K_{i+1}=\{1\}$, it follows that $v \in V(\Lambda_n)$ and for every $w\in K(l)$, $w\neq v$, there is $i\in [n]$ such that $w_i=1$ while $v_i=0$. Thus, $K(l)\cap V(\Lambda_n)=\{v\}$.

%\bigskip

%\noindent
%{\it Wydzia{\l} Matematyki, Informatyki i Ekonometrii, Uniwersytet Zielonog\'orski}\\
%{\it ul. Z. Szafrana 4a, 65-516 Zielona G\'ora, Poland}\\
%{\it E-mail: A.Kisielewicz@wmie.uz.zgora.pl}\\

\end{document}